\documentclass[12pt,twoside]{amsart} 
 \title{On weak Fano varieties with log canonical singularities}
\author{Yoshinori  Gongyo}
\address{Graduate School of Mathematical Sciences, 
the University of Tokyo, 3-8-1 Komaba, Meguro-ku, Tokyo 153-8914, Japan.}
\email{gongyo@ms.u-tokyo.ac.jp}
\date{\today}

\newcommand{\Supp}[0]{{\operatorname{Supp}}}

\usepackage{latexsym}
\usepackage{amsmath}
\usepackage{amssymb}
\usepackage{amsthm}
\usepackage{amscd}
\usepackage{enumerate}
\usepackage{amssymb} 
\usepackage{mathrsfs}          
\usepackage[all]{xy}

\usepackage{latexsym}

\newtheorem{thm}{Theorem}[section]

\newtheorem{prop}[thm]{Proposition}

\newtheorem{lem}[thm]{Lemma}

\newtheorem{rem}[thm]{Remark}

\newtheorem{cor}[thm]{Corollary}

\newtheorem{conj}[thm]{Conjecture}

\newtheorem{cl}[thm]{Claim}
 
\theoremstyle{definition}

\newtheorem{defi}[thm]{Definition}

\newtheorem{eg}[thm]{Example}

\newtheorem{defandlem}[thm]{Definition and Lemma}

\newtheorem{basic}[thm]{Basic construction}


\begin{document}
\bibliographystyle{amsalpha+}

 \maketitle
  \begin{abstract}We prove that the anti-canonical divisors of weak Fano $3$-folds with log canonical singularities are semi-ample. Moreover, we consider semi-ampleness of the anti-log canonical divisor of any weak log Fano pair with log canonical singularities. We show semi-ampleness dose not hold in general by constructing several examples. Based on those examples, we propose sufficient conditions which seem to be the best possible and we prove semi-ampleness under such conditions. In particular we derive semi-ampleness of the anti-canonical divisors of log canonical weak Fano $4$-folds whose lc centers are at most $1$-dimensional. We also investigate the Kleiman-Mori cones of weak log Fano pairs with log canonical singularities. 
 \end{abstract}

\section{Introduction} Throughout this paper, we work over $\mathbb{C}$, the complex number field. We start by some basic definitions.
\begin{defi}Let $X$ be a normal projective variety and $\Delta$ an effective $\mathbb{Q}$-Weil divisor on $X$. We say that $(X,\Delta)$ is a \emph{weak log Fano pair} if $-(K_X+\Delta)$ is nef and big. If $\Delta=0$, then we simply say that $X$ is a \emph{weak Fano variety}.
\end{defi}
\begin{defi}
Let $X$ be a normal variety and $\Delta$ an effective $\mathbb{Q}$-Weil divisor on $X$ such that $K_X+\Delta$ is a $\mathbb{Q}$-Cartier divisor. Let $\varphi:Y\rightarrow X$ be a log resolution of $(X,\Delta)$. We set $$K_Y=\varphi^*(K_X+\Delta)+\sum a_iE_i,$$ where $E_i$ is a prime divisor.
The pair $(X,\Delta)$ is called 
\begin{itemize}
\item[(a)] \emph{kawamata log terminal} $($\emph{klt}, for short$)$ if $a_i > -1$ for all $i$, or
\item[(b)]\emph{log canonical} $($\emph{lc}, for short$)$ if $a_i \geq -1$ for all $i$.
\end{itemize}
\end{defi}

\begin{defi}[Lc center]\label{lccenter} Let $(X, \Delta)$ be an lc pair.
We call that $C\subset X$ is an \emph{lc center} of $(X,\Delta)$ if there exists a log resolution $\varphi: Y\to X$ such that $\varphi(E) =C$ for some prime divisor $E$ on $Y$ with $a(E, X, \Delta)=-1$.
\end{defi}

There are questions whether the following fundamental properties hold or not for a log canonical weak log Fano pair $(X,\Delta)$  $($cf.\ \cite[2.6. Remark-Corollary]{S}, \cite[11.1]{P}$)$:
\begin{itemize}
\item[(i)]
Semi-ampleness of $-(K_X+\Delta)$.
\item[(ii)]
Existence of $\mathbb{Q}$-complements, i.e., existence of an effective $\mathbb{Q}$-divisor $D$ such that $K_X+\Delta+D \sim_{\mathbb{Q}} 0$ and $(X,\Delta+D)$ is lc.
\item[(iii)]Rational polyhedrality of the Kleiman-Mori cone $\overline{NE}(X)$.
\end{itemize}
It is easy to see that $($i$)$ implies $($ii$)$. In the case where $(X,\Delta)$ is a klt pair,  the above three properties hold by the Kawamata-Shokurov base point free theorem and the cone theorem (cf.\ \cite{KMM}, \cite{KoM}). Shokurov proved that these three properties hold for surfaces (cf.\ \cite[2.5. Proposition]{S}).  

Among other things, we prove the following:
\begin{thm}[={\rm{Corollaries}} \ref{Cor2} and \ref{Cor4}]\label{main}
Let $X$ be a weak Fano $3$-fold with log canonical singularities. 
Then $-K_X$ is semi-ample and $\overline{NE}(X)$ is a rational polyhedral cone.
\end{thm}
\begin{thm}[={\rm{Corollary}} \ref{Cor2-1} and {\rm{Theorem}} \ref{Thm3}]\label{main2}
Let $X$ be a weak Fano $4$-fold with log canonical singularities. Suppose  that any lc center of $X$ is at most $1$-dimensional.  
Then $-K_X$ is semi-ample and $\overline{NE}(X)$ is a rational polyhedral cone.
\end{thm}
On the other hand, the above three properties do not hold for $d$-dimensional log canonical weak \emph{log} Fano pairs in general, where $d\geq 3$. Indeed, we give the following examples of plt weak log Fano pairs whose anti-log canonical divisors are not semi-ample in Section \ref{counterexample} (in particular, such examples of $3$-dimensional weak log Fano plt pairs show the main result of \cite{Kar1} does not hold). It is well known that there exists a $(d-1)$-dimensional smooth projective variety $S$ such that $-K_S$ is nef and is not semi-ample (e.g. When $d=3$, we take a very general $9$-points blow up of $\mathbb{P}^2$ as $S$). Let $X_0$ be the cone over $S$ with respect to some projectively normal embedding $S\subset \mathbb{P}^N$. We take the blow-up $X$ of $X_0$ at its vertex. Let $E$ be the exceptional divisor of the blow-up. Then the pair $(X,E)$ is a weak log Fano plt pair such that $-(K_X+E)$ is not semi-ample. Moreover we give an example of a log canonical weak log Fano pair without $\mathbb{Q}$-complements and an example whose Kleiman-Mori cone is not polyhedral.

We now outline the proof of semi-ampleness of $-K_X$ as in Theorem $\ref{main}$. First, we take a birational morphism $\varphi:Y \rightarrow X$ such that $\varphi^*(K_X)=K_Y+S$, $(Y,S)$ is dlt and $S$ is reduced. We set $C:=\varphi(S)$, which is the union of lc centers of $X$. By an argument in the proof of the Kawamata-Shokurov base point free theorem (Lemma \ref{lem3}), it is sufficient to prove that $-(K_Y+S)|_{S}$ is semi-ample. Moreover we have only to prove that $-K_X|_{C}$ is semi-ample by the formula $K_X|_{C}=(\varphi|_{S})^*((K_Y+S)|_{S})$. 

It is not difficult to see semi-ampleness of the restriction of $-K_X$ on any lc center of $X$. The main difficulty is how to extend semi-ampleness to $C$ from each $1$-dimensional irreducible component $C_i$ of $C$ since
the configuration of $C_i$'s may be complicated.
The key to overcome this difficulty is the abundance theorem for 
$2$-dimensional semi-divisorial log terminal pairs (\cite{AFKM}). We decompose $C=C'\cup C''$, where
$$\Sigma:=\{i|\ -K_X|_{C_i}\equiv 0\}, \ C':=\bigcup_{i \in \Sigma}C_i,\ \text{and}\ \ C'':=\bigcup_{i \not\in \Sigma}C_i.$$
Let $S'$ be the union of the irreducible components of $S$
whose image on $X$ is contained in $C'$.
We define the boundary $\mathrm{Diff}_{S'} (S)$ on $S'$ by the
formula $K_Y+S|_{S'}=K_{S'}+\mathrm{Diff}_{S'} (S)$.
The pair $(S',\mathrm{Diff}_{S'} (S))$ is known to be
semi-divisorial log terminal pair (sdlt, for short).
Applying the abundance theorem to the pair
$(S',\mathrm{Diff}_{S'} (S))$, we see that
$K_{S'}+\mathrm{Diff}_{S'} (S)$ is $\mathbb{Q}$-linearly trivial, namely,
there is a non-zero integer $m_1$ such that 
$-m_1 (K_Y+S)|_{S'}=-m_1(K_{S'}+\mathrm{Diff}_{S'} (S))\sim 0$.
This shows that $-m_1 K_X|_{C'}\sim 0$.
On the other hand, since $-K_X|_{C''}$ is ample, we can take enough
sections of $H^0(C'',-m_2 K_X|_{C''})$ for a sufficiently large and divisible $m_2$ (Lemma \ref{lem2-1}). Thus, we can find enough sections of
$H^0(C,-{m} K_X|_{C})$ for a sufficiently large and divisible $m$,
and can conclude that  $-K_X|_C$ is semi-ample.

To generalize this theorem to higher dimensional weak log Fano pairs, let us recall the following conjectures:
 \begin{conj}[Abundance conjecture in a special case]\label {Conj2}
Let $(X, \Delta)$ be a $d$-dimensional projective sdlt pair whose $K_X +\Delta$ is numerically trivial. Then $K_X+\Delta$ is $\mathbb{Q}$-linearly trivial, i.e., there exists an $n \in \mathbb{N}$ such that $n(K_X +\Delta) \sim 0$. 
\end{conj}

The abundance conjecture is one of the most famous conjecture in the minimal model program. This conjecture is true when $d\leq3$ by the works of Fujita, Kawamata, Miyaoka, Abramovich, Fong, Koll\'ar, $\mathrm{M^{c}}$Kernan, Keel, Matsuki, and Fujino. 

By the same way as in the $3$-dimensional case, we see the following theorem:
\begin{thm}[=Theorem\ \ref{Thm1}]\label{main2} Assume that $\mathrm{Conjecture}\ \ref{Conj2}$ in dimension $d-1$ holds. Let $(X, \Delta)$ be a $d$-dimensional log canonical weak log Fano pair. Suppose that $M(X, \Delta)\leq 1$, where 
$$M(X, \Delta):=\mathrm{max} \{ \dim P|\ \text{$P$ is an lc center of $(X,\Delta)$} \}.$$
Then $-(K_X+\Delta)$ is semi-ample.
\end{thm}
Indeed, semi-ampleness of $-K_X$ as in Theorem \ref{main} is derived from the above theorem since the singular locus of any normal $3$-fold is at most $1$-dimensional and Conjecture \ref{Conj2} for surfaces holds (\cite{AFKM}). We also derive semi-ampleness of weak Fano $4$-folds such that $M(X,0) \leq 1$ because Conjecture \ref{Conj2} for $3$-folds holds (\cite{F3}). We remark that by $\mathrm{Examples}\ \ref{eg1}$ and $\ref{eg1-2}$, this condition for the dimension of lc centers is the best possible.

In Section \ref{cone}, by the cone theorem for normal varieties by Ambro and Fujino (cf.\ Theorem\ \ref{thm144}), we derive the following:
\begin{thm}[=Theorem\ \ref{Thm3}] Let $(X, \Delta)$ be a $d$-dimensional log canonical weak log Fano pair. Suppose that $M(X, \Delta)\leq 1$. Then $\overline {NE}(X)$ is a rational polyhedral cone.
\end{thm}
Note that rational polyhedrality of $\overline{NE}(X)$ as in Theorem \ref{main} is a corollary of the above theorem. In Example \ref{eg3}, we also see that the Kleiman-Mori cone is not rational polyhedral in general when $M(X,\Delta)\geq 2$.

This paper is based on the minimal model theory for log canonical pairs developed by Ambro and Fujino (\cite{Am1}, \cite{Am4}, \cite{Am3}, \cite{F7}, \cite{F6}, \cite{F2}).

We will make use of the standard notation and definitions as in \cite{KoM}.

\subsection*{Acknowledgment}
The author wishes to express his deep gratitude to his supervisor Prof.\ Hiromichi Takagi for various comments and many suggestions and to Prof.\ Osamu Fujino for very helpful and essential advices. In particular, Prof.\ Fujino allows him to carry the proof of Theorem \ref{dltblowup} in this paper. He wishes to thank Prof.\ Yujiro Kawamata for warm encouragement and valuable comments, Prof.\ Shunsuke Takagi for informations about semi-normality.
He is indebted to Dr.\ Shinnosuke Okawa, Dr.\ Taro Sano, Dr.\ Kiwamu Watanabe and Dr.\ Katsuhisa Furukawa. He also would like to thank Dr.\ I. V. Karzhemanov for answering many questions. 

\section{Preliminaries and Lemmas}\label{Prelimi}
In this section, we introduce notation and some lemmas for the proof of Theorem \ref{main2} (=$\mathrm{Theorem}$ \ref{Thm1}). 

\begin{defi} For a $\mathbb Q$-Weil divisor 
$D=\sum _{j=1}^r d_j D_j$ such that 
$D_j$ is a prime divisor for 
every $j$ and $D_i\ne D_j$ for $i\ne j$, we define 
the {\em{round-up}} $\ulcorner D\urcorner =\sum _{j=1}^{r} 
\ulcorner d_j\urcorner D_j$ 
(resp.~the {\em{round-down}} $\llcorner D\lrcorner 
=\sum _{j=1}^{r} \llcorner d_j \lrcorner D_j$), 
where for every real number $x$, 
$\ulcorner x\urcorner$ (resp.~$\llcorner x\lrcorner$) is 
the integer defined by $x\leq \ulcorner x\urcorner <x+1$ 
(resp.~$x-1<\llcorner x\lrcorner \leq x$). 
The {\em{fractional part}} $\{D\}$ 
of $D$ denotes $D-\llcorner D\lrcorner$. 
We define 
\begin{align*}&
D^{=1}=\sum _{d_j=1}D_j, \ \ D^{\leq 1}=\sum_{d_j\leq 1}d_j D_j, \\ &
D^{<1}=\sum_{d_j< 1}d_j D_j, \ \ \text{and}\ \ \ 
D^{>1}=\sum_{d_j>1}d_j D_j. 
\end{align*}
We call $D$ a {\em{boundary}} 
$\mathbb Q$-divisor if 
$0\leq d_j\leq 1$ 
for every $j$. 
\end{defi}

\begin{defi}[Stratum] Let $(X, \Delta)$ be an lc pair. 
A {\em{stratum}} of $(X, \Delta)$ denotes 
$X$ itself or an lc center of $(X, \Delta)$. 
\end{defi}

The following theorem is very important as a generalization of vanishing theorems (cf.\ \cite[Theorem 3.1]{Am4}, \cite[Theorem 2.2]{F7}, \cite[Theorem 2.38]{F6}, \cite[Theorem 6.3]{F2}). 

\begin{thm}[Torsion-freeness theorem]\label{thm53}Let 
$Y$ be a smooth variety and 
$B$ a boundary $\mathbb R$-divisor such that 
$\mathrm{Supp} B$ is simple normal crossing. 
Let $f:Y\to X$ be a projective morphism and $L$ a Cartier 
divisor on $Y$ such that 
$H\sim _{\mathbb R}L-(K_Y+B)$ is $f$-semi-ample. Then every associated prime of $R^qf_*\mathcal O_Y(L)$ is the generic point of 
the $f$-image of some stratum of $(Y, B)$  for any non-negative integer $q$. 
\end{thm}

The following theorem is proved by Fujino (\cite[Theorem 10.5]{F2}). We include the proof for the reader's convenience.
\begin{thm}\label{dltblowup}
Let $X$ be a normal quasi-projective variety and 
$\Delta$ an effective $\mathbb Q$-divisor on $X$ such 
that $K_X+\Delta$ is $\mathbb Q$-Cartier. Suppose that $(X,\Delta)$ is lc.
Then there exists a projective birational 
morphism $\varphi:Y\to X$ from a normal quasi-projective 
variety with the following properties: 
\begin{itemize}
\item[(i)] $Y$ is $\mathbb Q$-factorial, 
\item[(ii)] $a(E, X, \Delta)= -1$ for every  
$\varphi$-exceptional divisor $E$ on $Y$, 
\item[(iii)] for $$
\Gamma=\varphi^{-1}_*\Delta+\sum _{E: {\text{$\varphi$-exceptional}}}E, 
$$ it holds that  $(Y, \Gamma)$ is dlt and $K_Y+\Gamma=\varphi^*(K_X+\Delta)$, and
\item[(iv)] Let $\{C_i\}$ be any set of lc centers of $(X,\Delta)$. Let $W=\bigcup C_i$ with a reduced structure and $S$ the union of the irreducible components of $\llcorner\Gamma \lrcorner$ which are mapped into $W$ by $\varphi$. Then $(\varphi|_{S})_*\mathcal{O}_S \simeq \mathcal{O}_W$.
\end{itemize}

\end{thm}
\begin{proof}
Let $\pi: V\rightarrow X$ be a resolution such that
\begin{itemize}
\item[(1)] $\pi^{-1}(C)$ is a simple normal crossing divisor on $V$ for every lc center $C$ of $(X,\Delta)$, and
\item[(2)] $\pi_*^{-1}\Delta \cup \mathrm{Exc}(\pi) \cup \pi^{-1}( \mathrm{Nklt}(X, \Delta))$ has a simple normal crossing support, where $\mathrm{Exc}(\pi)$ is the exceptional set of $\pi$ and $\mathrm{Nklt}(X, \Delta)$ is the union of lc centers of $(X, \Delta)$.
\end{itemize}
By Hironaka's resolution theorem, we can assume that $\pi$ is a composite of blow-ups with centers of codimension at least two. Then there exists an effective $\pi$-exceptional Cartier divisor $B$  on $V$ such that $-B$ is $\pi$-ample. We put 
$$
F=\sum _{\begin{matrix}\scriptstyle{a(E, X, \Delta)>-1, }\\ 
\scriptstyle {E: {\text{$\pi$-exceptional}}}\end{matrix}}E
\ \text{and}\ G=\sum _{a(E, X, \Delta)= -1}E. 
$$
Let $H$ be a sufficiently ample Cartier divisor on $X$ such that $-B+\pi^*H$ is ample. We choose $0<\varepsilon \ll 1$ such that $\varepsilon G-B +\pi^*(H)$ is ample. Since $-B +\pi^*(H)$ and $\varepsilon G-B +\pi^*(H)$ are ample, we can take effective $\mathbb{Q}$-divisors $H_1$ and $H_2$ on $V$ with small coefficients such that $G+F +\pi^{-1}_{*}\Delta +H_1+H_2$ has a simple normal crossing support and that $-B +\pi^*H \sim_{\mathbb{Q}} H_1$, $\varepsilon G-B +\pi^*(H) \sim_{\mathbb{Q}} H_2$. We take $0<\nu, \mu \ll 1$ such that every divisor in $F$ has a negative coefficient in 
$$ M:=\Gamma_{V} - G- (1-\nu) F- \pi^{-1}_*\Delta^{<1} +\mu B,$$
where $\Gamma_V$ is a $\mathbb{Q}$-divisor on $V$ such that $K_V+\Gamma_V = \pi^*(K_X+\Delta)$. Now we construct a log minimal model of $(V, G +(1-\nu)F + \pi^{-1}_*\Delta^{<1} + \mu H_1)$ over $X$. Since 
$$G+(1-\nu)F +\mu H_1 \sim_{\mathbb{Q}}(1-\varepsilon \mu)G+(1-\nu)F +\mu H_2,$$
it is sufficient to construct a log minimial model of $(V, (1-\varepsilon\mu)G +(1-\nu)F + \pi^{-1}_*\Delta^{<1} + \mu H_2)$ over $X$. Because $(V, (1-\varepsilon\mu)G +(1-\nu)F + \pi^{-1}_*\Delta^{<1} + \mu H_2)$ is klt, we can get a log minimal model $\varphi: Y \to X $ of $(V, (1-\varepsilon\mu)G +(1-\nu)F + \pi^{-1}_*\Delta^{<1} + \mu H_2)$ over $X$ by \cite[Theorem 1.2]{BCHM}.

We show this $Y$ satisfies the conditions of the theorem. For any divisor $D$ on $V$ (appearing above), let $D'$ denote its strict transform on $Y$. We see the following claim:
\begin{cl}\label{f-cl2}
$F'=0$.
\end{cl}
\begin{proof}[Proof of Claim \ref{f-cl2}]
By the above construction,
$$N:=K_Y+G' +(1-\nu)F' + \varphi^{-1}_*\Delta^{<1} + \mu H_1'
$$
is $\varphi$-nef. Then 
\begin{eqnarray*}-M' &\sim_{\mathbb{Q}, \varphi}& N-(K_Y+\Gamma_Y)
\end{eqnarray*}
since $(\pi^*H)'=\varphi^*H$, hence it is $\varphi$-nef. Since $\varphi_{*}M'=0$, we see that $M'$ is effective by the negativity lemma (cf.\ \cite[Lemma 3.39]{KoM}). Since every divisor in $F$ has a negative coefficient in $M$, $F$ is contracted on $Y$. We finish the proof of Claim \ref{f-cl2}.
\end{proof}
From Claim \ref{f-cl2}, the discrepancy of every $\varphi$-exceptional divisor is equal to $-1$. We see that $Y$ satisfies the condition (ii). By the above construction, $(Y,\Gamma)$ is a $\mathbb{Q}$-factorial dlt pair since so is $(Y, G' + \varphi^{-1}_*\Delta^{<1} + \mu H_1)$. We see the condition (i). Because the support of $K_Y+\Gamma-\varphi^*(K_X+\Delta)$ coincide with $F'$, we see the condition (iii).

Now, we show that $Y$ and $\varphi$ satisfy the condition (iv). Since we get $Y$ by the log minimal model program over $X$ with scaling of some effective divisor with respect to $K_V+ G +(1-\nu)F + \pi^{-1}_*\Delta^{<1} + \mu H_1$ (cf.\ \cite{BCHM}), we see that the rational map $f: V \dashrightarrow Y$ is a composition of $(K_V+ G +(1-\nu)F + \pi^{-1}_*\Delta^{<1} + \mu H_1)$-negative divisorial contractions and log flips. Let $\Sigma$ be an lc center of $(Y, \Gamma)$. Then it is also an lc center of $(Y, \Gamma + \mu H_1')$. By the negativity lemma, 
$f: V \dashrightarrow Y$ is an isomorphism around the generic point of $\Sigma$. Therefore, if 
$\varphi(\Sigma) \subseteq W$, then $\Sigma \subseteq S$ by the conditions (1) and (2) for $\pi:V \to X$. This means that no lc centers of $(Y, \Gamma-S)$ are mapped into $W$ by $\varphi$. Let $g:Z \to Y$ be a resolution such that 
\begin{itemize}
\item[(a)] $\Supp\ \Gamma_Z$ is a simple normal crossing divisor, where $\Gamma_Z$ is defined by $K_Z+\Gamma_Z=g^*(K_Y+\Gamma),$ and
\item[(b)] $g$ is an isomorphism over the generic point of any lc center of $(Y, \Gamma)$.
\end{itemize}
Let $S_Z$ be the strict transform of $S$ on $Z$. We consider the following short exact sequence
\begin{align*}\label{hhh}
0\to \mathcal O_Z(\ulcorner -(\Gamma^{<1}_Z)\urcorner -S_Z)&\to \mathcal O_Z(\ulcorner 
-(\Gamma^{<1}_Z)\urcorner ) \tag{$\ast$}  \\ &\to \mathcal O_{S_Z}(\ulcorner -(\Gamma^{<1}_Z)\urcorner)\to 0. 
\end{align*}
We note that 
$$\ulcorner -(\Gamma^{<1}_Z)\urcorner -S_Z-(K_Z+\{\Gamma_Z\}+\Gamma_Z^{=1}-S_Z) \sim_{\mathbb{Q}} -h^*(K_X+\Delta), 
$$
where $h= \varphi \circ g$. Then we obtain 
\begin{align*}
0&\to h_*\mathcal O_Z(\ulcorner -(\Gamma^{<1}_Z)\urcorner -S_Z)\to h_*\mathcal O_Z(\ulcorner 
-(\Gamma^{<1}_Z)\urcorner )\to h_*\mathcal O_{S_Z}(\ulcorner -(\Gamma^{<1}_Z)\urcorner)\\ &\overset {\delta}\to 
R^1h_*\mathcal O_Z(\ulcorner -(\Gamma^{<1}_Z)\urcorner-S_Z)\to \cdots. 
\end{align*}
We claim the following:
\begin{cl}\label{f-cl1} $\delta$ is a zero map.
\end{cl}
\begin{proof}[Proof of Claim \ref{f-cl1}] Let $\Sigma$ be an lc center of $(Z, \{\Gamma_Z\}+\Gamma^{=1}_Z-S_Z)$. Then $\Sigma$ is some intersection of components of $\Gamma^{=1}_Z-S_Z$. By the conditions (a) and (b), $\Gamma^{=1}_Z-S_Z$ is the strict transform of $\llcorner \Gamma \lrcorner -S$. By this, the image of $\Sigma$ by $g$ is some intersection of components of $\llcorner \Gamma \lrcorner -S$. In particular, $g(\Sigma)$ is an lc center of $(Y, \Gamma-S)$. Thus no lc centers of $(Z, \{\Gamma_Z\}+\Gamma^{=1}_Z-S_Z)$ are mapped into $W$ by $h$.
Assume by contradiction that $\delta$ is not zero. Then there exists a section $s \in H^0(U, h_*\mathcal O_{S_Z}(\ulcorner -(\Gamma^{<1}_Z)\urcorner))$ for some non-empty open set $U\subseteq X$ such that $\delta(s)\not=0$. Since $\Supp\ \delta(s) \not= \emptyset$, we can take an associated prime $x \in \Supp\ \delta(s)$. We see that $x \in W$ since $\Supp (h_*\mathcal O_{S_Z}(\ulcorner -(\Gamma^{<1}_Z)\urcorner))$ is contained in $W$. By Theorem \ref{thm53}, $x$ is the generic point of the $h$-image of some stratum of $(Z, \{\Gamma_Z\}+\Gamma^{=1}_Z-S_Z)$. Since $h$ is a birational morphism, $x$ is the generic point of the $h$-image of some lc center of $(Z, \{\Gamma_Z\}+\Gamma^{=1}_Z-S_Z)$. Because no lc centers of $(Z, \{\Gamma_Z\}+\Gamma^{=1}_Z-S_Z)$ are mapped into $W$ by $h$, it holds that $x \not \in W$. But this contradicts the way of taking $x$.
 
\end{proof}
Thus, we obtain
$$0\to \mathcal{I} _W \to \mathcal{O}_X \to h_*\mathcal{O}_{S_Z}(\ulcorner -(\Gamma^{<1}_Z)\urcorner)\to 0,
$$
where $\mathcal{I} _W$ is the defining ideal sheaf of $W$ since $\ulcorner -(\Gamma^{<1}_Z)\urcorner$ is effective and $h$-exceptional.
This implies that $\mathcal{O}_W\simeq  h_*\mathcal{O}_{S_Z}(\ulcorner -(\Gamma^{<1}_Z)\urcorner)$. By applying $g_*$ to (\ref{hhh}), we obtain 
$$0\to \mathcal{I} _S \to \mathcal{O}_Y \to g_*\mathcal{O}_{S_Z}(\ulcorner -(\Gamma^{<1}_Z)\urcorner)\to 0,
$$
where $\mathcal{I} _S$ is the defining ideal sheaf of $S$ since $\ulcorner -(\Gamma^{<1}_Z)\urcorner$ is effective and $g$-exceptional. We note that
$$R^1g_*\mathcal{O}_Z(\ulcorner -(\Gamma^{<1}_Z)\urcorner-S_Z)=0
$$
by Theorem \ref{thm53} since $g$ is an isomorphism at the generic point of any stratum of $(Z, \{\Gamma_Z\}+\Gamma_Z^{=1}-S_Z)$.
Thus, $\mathcal{O}_W\simeq h_*\mathcal{O}_{S_Z}(\ulcorner -(\Gamma^{<1}_Z)\urcorner) \simeq \varphi_*g_*\mathcal{O}_{S_Z}(\ulcorner -(\Gamma^{<1}_Z)\urcorner) \simeq \varphi_*\mathcal{O}_S$. We finish the proof of  Theorem \ref{dltblowup}.
\end{proof}

\begin{defi} Let $X$ be a normal variety and $D$ a $\mathbb{Q}$-Weil divisor. We define that $$R(X, D)=\bigoplus_{m=0}^{\infty} H^0(X, \llcorner mD \lrcorner).$$
 
\end{defi}

\begin{defi}[semi-divisorial log terminal, cf.\ \cite{F3}]Let $X$ be a reduced $S_2$-scheme. We assume that it is pure $d$-dimensional and is normal crossing in codimension $1$. Let $\Delta$ be an effective $\mathbb{Q}$-Weil divisor on $X$ such that $K_X+\Delta$ is $\mathbb{Q}$-Cartier. 

Let $X=\bigcup X_i$ be the decomposition into irreducible components, and $\nu:X':=\coprod X'_i \rightarrow X=\bigcup X_i$ the normalization. Define the $\mathbb{Q}$-divisor $\Theta$ on $X'$ by $K_{X'}+\Theta:=\nu^*(K_X+\Delta)$ and set $\Theta_i:=\Theta|_{X_i'}$. 

We say that $(X,\Delta)$ is \emph{semi-divisorial log terminal} (for short, \emph{sdlt}) if $X_i$ is normal, that is, $X_i'$ is isomorphic to $X_i$, and $(X_i',\Theta_i)$ is dlt for every $i$. 
 \end{defi}

\begin{defandlem}[Different, cf.\ \cite{C}]\label{different} Let $(Y,\Gamma)$ be a dlt pair and $S$ a union of some components of $\llcorner\Gamma\lrcorner$. Then there exists an effective $\mathbb{Q}$-divisor $\mathrm{Diff}_{S}(\Gamma)$ on $S$ such that $(K_Y+\Gamma)|_{S}\sim_{\mathbb{Q}} K_S+\mathrm{Diff}_S(\Gamma)$. The effective $\mathbb{Q}$-divisor $\mathrm{Diff}_{S}(\Gamma)$ is called the \emph{different} of $\Gamma$. Moreover it holds that $(S, \mathrm{Diff}_{S}(\Gamma))$ is sdlt.
\end{defandlem}

The following proposition is \cite [Proposition 2]{fk2} 
(for the proof, see \cite [Proof of Theorem 3]{fk1} and \cite 
[Lemma 3]{Ka3}). 

\begin{prop}\label{prop}
Let $(X,\Delta)$ be a proper dlt pair and 
$L$ a nef Cartier divisor such that $aL-(K_X+\Delta)$ is nef 
and big for some $a\in {\mathbb N}$. 
If $\mathrm{Bs} |mL| \cap \llcorner \Delta\lrcorner=\emptyset$ for 
every $m\gg 0$, then $|mL|$ is base point free for every $m\gg 0$, 
where $\mathrm{Bs} |mL|$ is the base locus of $|mL|$. 
\end{prop}

By this proposition, we derive the following lemma:

\begin{lem}\label{lem3} Let $(Y,\Gamma)$ be a $\mathbb{Q}$-factorial weak log Fano dlt pair. Suppose that $-(K_{S} +\Gamma_{S})$ is semi-ample, where $S:=\llcorner\Gamma\lrcorner$ and $\Gamma_{S}:=\mathrm{Diff}_{S}(\Gamma)$. Then $-(K_Y + \Gamma)$ is semi-ample.
\end{lem}
\begin{proof}We consider the exact sequence
\begin{eqnarray}
\nonumber 0 \to \mathcal{O}_{Y}(-m(K_{Y}+\Gamma)-S)
\to \mathcal{O}_{Y}(-m(K_{Y}+\Gamma)) \to \\
\nonumber \to
\mathcal{O}_{S}(-m(K_{Y}+\Gamma)|_{S}) \to 0
\end{eqnarray}
for $m \gg 0$. By the Kawamata-Viehweg vanishing theorem (cf.\ \cite[Theorem 1-2-5.]{KMM}, \cite[Theorem 2.70]{KoM}), we have
\begin{eqnarray}
\nonumber H^{1}(Y,\mathcal{O}_{Y}(-m(K_{Y}+\Gamma)-S)) = \\
\nonumber = H^{1}(Y,\mathcal{O}_{Y}(K_{Y} +\Gamma-S -
(m+1)(K_{X}+\Gamma))) = \{0 \},
\end{eqnarray}
since the pair $(Y,\Gamma-S)$ is klt and $-(K_Y + \Gamma)$ is nef and big. Thus, we get the 
exact sequence
\begin{eqnarray}
\nonumber H^{0}(Y,\mathcal{O}_{Y}(-m(K_{Y}+\Gamma)) \to
H^{0}(S,\mathcal{O}_{S}(-m(K_{Y}+\Gamma)|_{S}))
\to 0.\nonumber
\end{eqnarray}
Therefore, we see that $\mathrm{Bs} |-m(K_Y+\Gamma)| \cap S =\emptyset$ for $m \gg 0$ since $-(K_S+\Delta_S)$ is semi-ample. Applying $\mathrm{Proposition}\ \ref{prop}$, we conclude that $-(K_Y + \Gamma)$ is semi-ample.
\end{proof}

\begin{defi}(cf.\ \cite[1.1. Definition]{GT}, \cite[Definition 7.1]{KoS})\label{s-int} Suppose that $R$ is a reduced excellent ring and $R\subseteq S$ is a reduced $R$-algebra which is finite as an $R$-module. We say that the extension $i:R\hookrightarrow S$ is \emph{subintegral} if one of the following equivalent conditions holds:
\begin{itemize}
\item[(a)] $(S\bigotimes_R k(\mathfrak{p}))_{\mathrm{red}}=k(\mathfrak{p})$ for all $\mathfrak{p} \in \mathrm{Spec}(R)$. 
\item[(b)] the induced map on the spectra is bijective and $i$ induces trivial residue field extensions.
\end{itemize}
\end{defi}

\begin{defi}\cite[Definition 7.2]{KoS}\label{semi-normal} Suppose that $R$ is a reduced excellent ring. We say that $R$ is \emph{semi-normal} if every subintegral extension $R\hookrightarrow S$ is an isomorphism. 

A scheme $X$ is called \emph{semi-normal at} $q \in X$ if the local ring at $q$ is semi-normal. 
If $X$ is semi-normal at every point, we say that $X$ is \emph{semi-normal}.
\end{defi} 
 
\begin{prop}\cite[5.3. Corollary]{GT}\label{ff} Let $(R,\mathfrak{m})$ be a local excellent ring. Then $R$ is semi-normal if and only if $\widehat{R}$ is semi-normal, where $\widehat{R}$ is $\mathfrak{m}$-adic completion of $R$.
\end{prop} 
 
\begin{prop}$($cf.\ \cite[7.2.2.1]{Ko}, \cite[Remark 7.6]{KoS}$)$\label{prop4} Let $C$ be a pure $1$-dimensional proper reduced scheme of finite type over $\mathbb{C}$, and $q \in C$ a closed point. Then $C$ is semi-normal at $q$ if and only if $\widehat{\mathcal{O}}_{C,q}$ satisfies that
 \begin{itemize}
\item[(i)] $\widehat{\mathcal{O}}_{C,q}\simeq \mathbb{C}[[X]]$, or 
\item[(ii)] $\widehat{\mathcal{O}}_{C,q}\simeq \mathbb{C}[[X_1,X_2,\cdots ,X_r]]/ \langle X_iX_j|1\leq i \not=j \leq r \rangle$ for some $r\geq 2$, i.e., $q \in C$ is isomorphic to the coordinate axies in $\mathbb{C}^r$ at the origin as a formal germs.  
 \end{itemize} 
\end{prop}

\begin{lem}\label{lem2-1} Let $C=C_1\cup C_2$ be a pure $1$-dimensional proper semi-normal reduced scheme of finite type over $\mathbb{C}$, where $C_1$ and $C_2$ are pure $1$-dimensional reduced closed subschemes.  Let $D$ be a $\mathbb{Q}$-Cartier divisor on $C$.
Suppose that  $D_1$ is $\mathbb{Q}$-linearly trivial and $D_2$ is ample, where $D_i:=D|_{C_i}$. 
Then $D$ is semi-ample.

\end{lem}
\begin{proof} Let $C_1 \cap C_2 =\{p_1, \dots, p_r\}$. We take $m \gg 0$ which satisfies the following:
\begin{itemize}
\item[(i)] $mD_1 \sim 0$,
\item[(ii)] $\mathcal{O}_{C_2}(mD_2) \otimes (\bigcap_{k \not = l}\mathfrak{m}_{p_k})$ is generated by global sections for all $l\in \{1,\dots, r\}$, and
\item[(iii)] $\mathcal{O}_{C_2}(mD_2) \otimes ( \bigcap_{k}\mathfrak{m}_{p_k})$ is generated by global sections,
\end{itemize}
where $\mathfrak{m}_{p_k}$ is the ideal sheaf of $p_k$ on $C_2$. We choose a nowhere vanishing
section $s\in H^0(C_1,mD_1)$.
By (ii), we can take a section $t_{l}\in H^0(C_2, mD_2)$
which does not vanish at $p_{l}$ but
vanishes at all the $p_{k}$ $(k\in \{1,\dots, r\}, k\not =l)$
for each $l\in \{1,\dots, r\}$.
By multiplying suitable nonzero constants to $t_{l}$,
we may assume that $t_{l}|_{p_{l}}=s|_{p_{l}}$. We set $t:=\sum_{l} t_l \in H^0(C_2, mD_2)$.
Since $C$ is semi-normal, $\mathrm{Proposition}\ \ref{prop4}$ implies that $\mathcal{O}_{C_1 \cap C_2} \simeq \bigoplus_{l=1}^{r}\mathbb{C}(p_l)$, where $\mathbb{C}(p_l)$ is the skyscraper sheaf $\mathbb{C}$ sitting at $p_l$, by computations on $\widehat{\mathcal{O}}_{C,p_l}$. Thus we get the following exact sequence:
$$0\to \mathcal{O}_C(mD) \to \mathcal{O}_{C_1}(mD_1)\oplus \mathcal{O}_{C_2}(mD_2)\to \bigoplus_{l=1}^{r}\mathbb{C}(p_l) \to 0,
$$
where the third arrow maps $(s', s'')$ to $((s'-s'')|_{p_1}, \dots ,(s'-s'')|_{p_r}).$
Hence $s$ and $t$ patch together and give a section $u$ of 
$H^0(C, mD)$.

Let $p$ be any point of $C$.
If $p\in C_1$, then $u$ does not vanish at $p$.
We may assume that $p\in C_2 \setminus C_1$. By (iii), we can take a section $t' \in H^0(C_2, mD_2)$ which does not vanish at $p$ but
vanishes at $p_{l}$ for all $l \in \{1,\dots, r\}$. The zero section $0 \in H^0(C_1, mC_1)$ and $t'$ patch together and give a section $u'$ of 
$H^0(C, mD)$. By construction, the section $u'$ does not vanish at $p$. We finish the proof of $\mathrm{Lemma}\ \ref{lem2-1}$.
\end{proof}

\section{On semi-ampleness for weak Fano varieties}\label{semi-ample} In this section, we prove Theorem \ref{main2} (=$\mathrm{Theorem}$ \ref{Thm1}). As a corollary, we see that the anti-canonical divisors of weak Fano 3-folds with log canonical singularities are semi-ample. Moreover we derive semi-ampleness of the anti-canonical divisors of log canonical weak Fano $4$-folds whose lc centers are at most $1$-dimensional.

\begin{thm}\label{Thm1} Assume that $\mathrm{Conjecture}\ \ref{Conj2}$ in dimension $d-1$ holds. Let $(X, \Delta)$ be a $d$-dimensional log canonical weak log Fano pair. Suppose that $M(X, \Delta)\leq 1$, where 
$$M(X, \Delta):=\mathrm{max} \{ \dim P|\ \text{$P$ is an lc center of $(X,\Delta)$} \}.$$
Then $-(K_X+\Delta)$ is semi-ample.
\end{thm}
 
\begin{proof}By Theorem\ \ref{dltblowup}, we take a birational morphism $\varphi :(Y,\Gamma) \rightarrow (X,\Delta)$ as in the theorem. We set $S:=\llcorner\Gamma\lrcorner$ and $C:=\varphi(S)$, where we consider the reduced scheme structures on $S$ and $C$. We have only to prove that $-(K_{S} +\Gamma_{S})=-(K_Y + \Gamma)|_{S}$ is semi-ample from {\rm{Lemma}} \ref{lem3}. By the formula $(K_Y+\Gamma)|_{S} \sim_{\mathbb{Q}} (\varphi|_{S})^*((K_X+\Delta)|_{C})$, it suffices to show that $-(K_X+\Delta)|_{C}$ is semi-ample. Arguing on each connected component of $C$, we may assume that $C$ is connected. Since $M(X, \Delta)\leq 1$, it holds that $\dim C\leq 1$. When $\dim C=0$, i.e., $C$ is a closed point, then $-(K_X+\Delta)|_{C}\sim_{\mathbb{Q}} 0$, in particular, is semi-ample.\\ 
When $\dim C=1$, $C$ is a pure $1$-dimensional semi-normal scheme by \cite[Theorem 1.1]{Am3} or \cite[Theorem 9.1]{F2}. Let $C=\bigcup_{i=1}^{r}C_i$, where $C_i$ is an irreducible component, and let $D:=-(K_X+\Delta)|_{C}$ and $D_i:=D|_{C_i}$. 
We set
$$\Sigma:=\{i|\ D_i\equiv 0\}, \ C':=\bigcup_{i \in \Sigma}C_i ,\ C'':=\bigcup_{i \not\in \Sigma}C_i.$$
Let $S'$ be the union of irreducible components of $S$ whose image by $\varphi$ is contained in $C'$. We see that $K_{S'}+\Gamma_{S'}\equiv 0$, where $\Gamma_{S'}:=\mathrm{Diff}_{S'}(\Gamma)$. Thus it holds that $K_{S'}+\Gamma_{S'}\sim_{\mathbb{Q}} 0$ by applying \rm{Conjecture} \ref{Conj2} to $(S',\Gamma_{S'})$. Since $(\varphi|_{S'})_{*}\mathcal{O}_{S'} \simeq \mathcal{O}_{C'}$ by the condition (iv) in Theorem\ \ref{dltblowup}, it holds that $D|_{C'} \sim_{\mathbb{Q}} 0$.
We see that $D|_{C''}$ is ample since the restriction of $D$ on any irreducible component of $C''$ is ample. By $\mathrm{Lemma}\ \ref{lem2-1}$, we see that $D=-(K_X+\Delta)|_{C}$ is semi-ample. We finish the proof of $\mathrm{Theorem}\ \ref{Thm1}$.
\end{proof}

\begin{cor}\label{ring} Assume that $\mathrm{Conjecture}\ \ref{Conj2}$ in dimension $d-1$ holds.\\
Let $(X, \Delta)$ be a $d$-dimensional log canonical weak log Fano pair. Suppose that $M(X, \Delta)\leq 1$. Then $R(X, -(K_X+\Delta))$ is a finitely generated algebra over $\mathbb{C}$.
\end{cor}

\rm{Conjecture} \ref{Conj2} holds for surfaces and $3$-folds by \cite{AFKM} and \cite{F3}. Thus we immediately obtain the following corollaries:

\begin{cor}\label {Cor2}
Let $(X, \Delta)$ be a $3$-dimensional log canonical weak log Fano pair. Suppose that $\llcorner\Delta\lrcorner=0$.
Then $-(K_X+\Delta)$ is semi-ample and $R(X, -(K_X+\Delta))$ is a finitely generated algebra over $\mathbb{C}$.
In particular, if $X$ is a weak Fano $3$-fold with log canonical singularities, then $-K_X$ is semi-ample and $R(X, -K_X)$ is a finitely generated algebra over $\mathbb{C}$.

\end{cor}

\begin{cor}\label {Cor2-1}
Let $(X, \Delta)$ be a $4$-dimensional log canonical weak log Fano pair. Suppose that $M(X,\Delta) \leq 1$.
Then $-(K_X+\Delta)$ is semi-ample and $R(X, -(K_X+\Delta))$ is a finitely generated algebra over $\mathbb{C}$.
In particular, if $X$ is a log canonical weak Fano $4$-fold whose lc centers are at most $1$-dimensional, then $-K_X$ is semi-ample and $R(X, -K_X)$ is a finitely generated algebra over $\mathbb{C}$.
\end{cor}

\begin{rem}\label{rem2}
When $M(X, \Delta)\geq 2$, $-(K_X+\Delta)$ is not semi-ample and $R(X, -(K_X+\Delta))$ is not a finitely generated algebra over $\mathbb{C}$, in general $(\mathrm{Examples}$ $\ref{eg1}$ and $\ref{eg1-2})$. 
\end{rem}

\begin{rem}\label{rem5} Based on $\mathrm{Theorem}\ \ref{Thm1}$, we expect the following statement:
\begin{quote}Let $(X, \Delta)$ be an lc pair and $D$ a nef Cartier divisor. Suppose there is a positive number $a$ such that $aD-(K_X+\Delta)$ is nef and big. If it holds that $M(X,\Delta)\leq 1$, then $D$ is semi-ample.
\end{quote}
However, there is a counterexample for this statement due to Zariski $($cf.\ \cite[Remark 3-1-2]{KMM}, \cite{Z}$)$.

\end{rem}

\section{On the Kleiman-Mori cone for weak Fano varieties}\label{cone} In this section, we introduce the cone theorem for normal varieties by Ambro and Fujino and prove polyhedrality of the Kleiman-Mori cone for a log canonical weak Fano variety whose lc centers are at most $1$-dimensional.
We use the notion of the scheme $\mathrm{Nlc}(X,\Delta)$, whose underlying space is the set of non-log canonical singularities. For the scheme structure on $\mathrm{Nlc}(X,\Delta)$, we refer \cite[Section 7]{F2} and \cite{F5} in detail.

\begin{defi}$($\cite[Definition 16.1]{F2}$)$\label{def141}
Let $X$ be a normal variety and $\Delta$ an effective 
$\mathbb Q$-divisor on $X$ such 
that $K_X+\Delta$ is $\mathbb Q$-Cartier. 
Let $\pi:X\to S$ be a projective morphism. 
We put\index{$\overline{NE}(X/S)_{-\infty}$} 
$$
\overline{NE}(X/S)_{\mathrm{Nlc}(X, \Delta)}=
\\ \mathrm{Im} (\overline {NE}(\mathrm{Nlc}(X, \Delta)/S)
\to \overline {NE}(X/S)). 
$$
\end{defi}

\begin{defi}$($\cite[Definition 16.2]{F2}$)$\label{def142}
An {\em{extremal face}}\index{extremal 
face} of $\overline {NE}(X/S)$ is a non-zero 
subcone $F\subset \overline {NE}(X/S)$ such that 
$z, z'\in F$ and $z+z'\in F$ implies that 
$z, z'\in F$. Equivalently, 
$F=\overline {NE}(X/S)\cap H^{\perp}$ for some 
$\pi$-nef $\mathbb R$-divisor $H$, which 
is called a {\em{supporting function}} of $F$.\index{supporting function} 
An {\em{extremal ray}}\index{extremal ray} 
is a one-dimensional 
extremal face. 
\begin{itemize}
\item[(1)] An extremal face $F$ is called 
{\em{$(K_X+\Delta)$-negative}}\index{$\omega$-negative} 
if $$F\cap \overline {NE}(X/S)_{K_X+\Delta\geq 0} 
=\{ 0\}. $$ 
\item[(2)] An extremal face $F$ is called {\em{rational}} 
if we can choose 
a $\pi$-nef $\mathbb Q$-divisor $H$ as a support 
function of $F$. 
\item[(3)] An extremal face $F$ is called {\em{relatively ample at 
$\mathrm{Nlc} (X, \Delta)$}}\index{relatively ample 
at infinity} if $$F\cap \overline {NE}(X/S)_{\mathrm{Nlc} (X, \Delta)}=\{0\}. 
$$
Equivalently, $H|_{\mathrm{Nlc} (X, \Delta)}$ is $\pi|_{\mathrm{Nlc} (X, \Delta)}$-ample for 
every supporting function $H$ of $F$. 
\item[(4)] An extremal face $F$ is called {\em{contractible 
at $\mathrm{Nlc} (X, \Delta)$}}\index{contractible at 
infinity} if it has a rational supporting function 
$H$ such that $H|_{\mathrm{Nlc} (X, \Delta)}$ is $\pi|_{\mathrm{Nlc} (X, \Delta)}$-semi-ample. 
\end{itemize}
\end{defi}

\begin{thm}{\rm{(Cone theorem for normal varieties, \cite[Theorem 5.10]{Am4}, \cite[Theorem 16.5]{F2})}}\label{thm144} 
Let $X$ be a normal variety, $\Delta$ an effective $\mathbb Q$-divisor 
on $X$ such that $K_X+\Delta$ is $\mathbb Q$-Cartier, and 
$\pi:X\to S$ a projective morphism. 
Then we have the following 
properties. 
\begin{itemize}
\item[$(1)$] $\overline {NE}(X/S)=\overline {NE}(X/S)_{K_X+\Delta\geq 0} 
+\overline {NE}(X/S)_{\mathrm{Nlc} (X, \Delta)}+\sum R_j$, 
where $R_j$'s are the $(K_X+\Delta)$-negative 
extremal rays of $\overline {NE}(X/S)$ that are 
rational and relatively ample at $\mathrm{Nlc} (X, \Delta)$. 
In particular, each $R_j$ is spanned by 
an integral curve $C_j$ on $X$ such that 
$\pi(C_j)$ is a point.  
\item[$(2)$] Let $H$ be a $\pi$-ample $\mathbb Q$-divisor 
on $X$. 
Then there are only finitely many $R_j$'s included in 
$(K_X+\Delta+H)_{<0}$. In particular, 
the $R_j$'s are discrete in the half-space 
$(K_X+\Delta)_{<0}$. 
\item[$(3)$] Let $F$ be a $(K_X+\Delta)$-negative extremal 
face of $\overline {NE}(X/S)$ that is 
relatively ample at $\mathrm{Nlc} (X, \Delta)$. 
Then $F$ is a rational face. 
In particular, $F$ is contractible at $\mathrm{Nlc} (X, \Delta)$. 
\end{itemize}
\end{thm}

By the above Theorem, we derive the following theorem:

\begin{thm}\label{Thm3} Let $(X, \Delta)$ be a $d$-dimensional log canonical weak log Fano pair. Suppose that $M(X, \Delta)\leq 1$. Then $\overline {NE}(X)$ is a rational polyhedral cone.
\end{thm}

\begin{proof}Since $-(K_X+\Delta)$ is nef and big, there exists an effective divisor $B$ satisfies the following: for any sufficiently small rational positive number $\varepsilon$, there exists a general $\mathbb{Q}$-ample divisor $A_{\varepsilon}$ such that $$-(K_X+\Delta)\sim_{\mathbb{Q}}\varepsilon B+A_{\varepsilon}.$$
We fix a sufficiently small rational positive number $\varepsilon$ and set $A:=A_{\varepsilon}$.
We also take a sufficiently small positive number $\delta$. Thus $\mathrm{Supp}(\mathrm{Nlc}(X,\Delta+\varepsilon B+\delta A))$ is contained in the union of lc centers of $(X,\Delta)$ and $-(K_X+\Delta+\varepsilon B +\delta A)$ is ample. By applying $\mathrm{Theorem}\ \ref{thm144}$ to $(X,\Delta+\varepsilon B +\delta A)$, We get
$$\overline {NE}(X)= \overline {NE}(X)_{\mathrm{Nlc} (X, \Delta+\varepsilon B+ \delta A)}+\sum_{j=1}^{m} R_j\ \text{for some $m$}.$$
Now we see that $\overline {NE}(X)_{\mathrm{Nlc} (X, \Delta+\varepsilon B+ \delta A)}$ is polyhedral since $\dim \mathrm{Nlc} (X, \Delta+\varepsilon B)\leq1$ by the assumption of $M(X,\Delta)\leq 1$. We finish the proof of $\mathrm{Theorem}\ \ref{Thm3}$.

\end{proof}

\begin{cor}\label {Cor4}
Let $X$ be a weak Fano $3$-fold with log canonical singularities. 
Then the cone $\overline {NE}(X)$ is rational polyhedral.\\
\end{cor}

\begin{rem}\label{rem3}
When $M(X, \Delta)\geq 2$, $\overline {NE}(X)$ is not polyhedral in general $(\mathrm{Example}\ \ref{eg3}).$
\end{rem}

\section{Examples}\label{counterexample}
In this section, we construct examples of log canonical weak log Fano pairs $(X,\Delta)$ such that $-(K_X+\Delta)$ is not semi-ample, $(X,\Delta)$ does not have $\mathbb{Q}$-complements, or $\overline{NE}(X)$ is not polyhedral.
 
\begin{basic}\label{basic}Let $S$ be a $(d-1)$-dimensional smooth projective variety such that $-K_S$ is nef and $S\subset\mathbb{P}^N$ some projectively normal embedding.
Let $X_0$ be the cone over $S$ and $\phi\colon X\to X_0$
the blow-up at the vertex.
Then 
the linear projection $X_0 \dashrightarrow S$ from the vertex
is decomposed as follows:
\begin{equation*}
\label{eq:line}
\xymatrix{ & X \ar[dl]_{\phi} \ar[dr]^{\pi}\\
 X_0 &  & S.}
\end{equation*}
This diagram is the restriction of 
the diagram for the projection $\mathbb{P}^{N+1}\dashrightarrow \mathbb{P}^N$:
\begin{equation*}
\label{eq:line}
\xymatrix{ & V:=\mathbb{P}_{\mathbb{P}^N}(\mathcal{O}_{\mathbb{P}^N}\oplus\mathcal{O}_{\mathbb{P}^N}(-1)) \ar[dl]_{\phi_0} 
\ar[dr]^{\pi_0}\\
 \mathbb{P}^{N+1} &  & \mathbb{P}^N.}
\end{equation*}
Moreover, the $\phi_0$-exceptional divisor is
the tautological divisor of $\mathcal{O}_{\mathbb{P}^N}\oplus\mathcal{O}_{\mathbb{P}^N}(-1)$.
Hence 
$X\simeq \mathbb{P}_S(\mathcal{O}_S\oplus\mathcal{O}_S(-H))$,
where $H$ is a hyperplane section on $S\subset\mathbb{P}^N$,
and 
the $\phi$-exceptional divisor $E$ is isomorphic to $S$ 
and is the tautological divisor of $\mathcal{O}_S\oplus\mathcal{O}_S(-H)$.
 
By the canonical bundle formula, 
it holds that
$$K_X= -2E+\pi^*(K_S-H),$$ thus we have
$$ -(K_X +E) =\pi^*(-K_S) + \pi^*H+E$$

We see $\pi^*H+E$ is nef and big since 
$\mathcal{O}_X(\pi^*(H) +E) \simeq \phi^*\mathcal{O}_{X_0}(1)$
and $\phi$ is birational. Hence $-(K_X + E)$ is nef and big since $\pi^*(-K_S)$ is nef. 
\end{basic}
The above construction is inspired by that of Hacon and $\mathrm{M^c}$Kernan in Lazi\'c's paper (cf.\ \cite[Theorem A.6]{L}).

In the following examples, $(X,E)$ is the plt weak log Fano pair given by the above construction.

\begin{eg}\label{eg1}This is an example of a $d$-dimensional plt weak log Fano pair such that the anti-log canonical divisors are not semi-ample, where $d\geq 3$.\\
There exists a variety $S$ such that $-K_S$ is nef and is not semi-ample (e.g. the surface obtained by blowing up $\mathbb{P}^2$ at very general $9$ points).
We see that $-(K_X + E)$ is not semi-ample since $-(K_X +E) |_E =-K_E$ is not semi-ample. In particular, $R(X,-(K_X +E))$ is not a finitely generated algebra over $\mathbb{C}$ by $-(K_X + \Delta)$ is nef and big.
\end{eg}
    
\begin{eg}\label{eg1-2}This is an example of a log canonical weak Fano variety such that the anti-canonical divisor is not semi-ample.\\
Let $T$ be a $k$-dimensional smooth projective variety whose $-K_{T}$ is nef and $A$ a $(d-k-1)$-dimensional smooth projective manifold with $K_A\sim_{\mathbb{Q}}0$, where $d$ and $k$ are integers satisfying $d-1 \geq k \geq 0$. We set $S=A \times T$. Let $p_{T}:S\rightarrow T$ be the canonical projection. We see that $K_S=p_{T}^{*}(K_T)$. Let $A_p$ be the fiber of $p_T$ at a point $p \in T$, and $\varphi:X\rightarrow Y$ the birational morphism with respect to $|\phi^* (\mathcal{O}_{X_0}(1))\otimes \pi^*p_T^*\mathcal{O}_T(H_T)|$, where $H_T$ is some very ample divisor on $T$. We claim the following:

\begin{cl}\label{cleg2} It holds that:
\begin{itemize}
\item[(i)] $Y$ is a projective variety with log canonical singularities.
\item[(ii)] $\mathrm{Exc}(\varphi)=E$ and any exceptional curve of $\varphi$ is contained in some $A_p$.
\item[(iii)] $\varphi^*K_Y=K_X+E$.
\item[(iv)] $\varphi(E) = T$ and $(\varphi|_{E})^*K_T=K_E$. 
\end{itemize}
\end{cl}
\begin{proof}[Proof of Claim\ \ref{cleg2}]
We see (ii) easily. Because  $-E|_{E}$ is ample, $E$ is not $\varphi$-numerical trivial. Set $\varphi^*K_Y=K_X+E+aE$ for some $a \in \mathbb{Q}$.  Since $K_X+E$ is $\varphi$-numerical trivial, we see $a=0$. Thus we obtain (iii). (i) follows from (iii).
By (iii), $\varphi(E)$ is an lc center. By $(\phi^* (\mathcal{O}_{X_0}(1))\otimes \pi^*p_T^*\mathcal{O}_T(H_T))|_{E}\simeq p_T^*\mathcal{O}_T(H_T)$, it holds that $\varphi|_{E}=p_T$. Thus (iv) follows.  
\end{proof}
If $-K_T$ is not semi-ample, then $-K_Y$ is not semi-ample and $k \geq 2$. 
Thus we see that $Y$ is a log canonical weak Fano variety with $M(Y,0)=k$ and $-K_Y$ is not semi-ample. In particular, $R(X,-K_X)$ is not a finitely generated algebra over $\mathbb{C}$ by $-K_X $ is nef and big (cf.\ \cite[Theorem 2.3.15]{Lf}).

\end{eg}

\begin{eg}\label{eg2} We construct an example of a weak log Fano plt pair without $\mathbb{Q}$-complements.\\
Let $S$ be the $\mathbb{P}^1$-bundle over an elliptic curve with respect to a non-split vector bundle of degree $0$ and  rank $2$. Then $-K_S$ is nef and $S$ does not have $\mathbb{Q}$-complements (cf.\ \cite[1.1. Example]{S}). Thus $(X,E)$ does not have $\mathbb{Q}$-complements by the adjunction formula $-(K_X +E) |_E =-K_E$.  

\end{eg}

\begin{eg}\label{eg3} We construct an example of a weak log Fano plt pair whose Kleiman-Mori cone is not polyhedral.
Let $S$ be the surface obtained by blowing up $\mathbb{P}^2$ at very general $9$ points. It is well known that $S$ has infinitely many $(-1)$-curves $\{C_i\}$.\\
Then we see that the Kleiman-Mori cone $\overline{NE}(X)$ is not polyhedral. Indeed, we have the following claim:
\begin{cl}\label{cleg}
$\mathbb{R}_{\geq0}[C_i] \subseteq \overline{NE}(X)$ is an extremal ray with $(K_X + E).C_i=-1$. Moreover, it holds that $\mathbb{R}_{\geq0}[C_i]\not=\mathbb{R}_{\geq0}[C_j]$ $(i\not=j)$.
\end{cl}
\begin{proof}[Proof of Claim\ \ref{cleg}]
We take a semi-ample line bundle $L_i$ on $S$ such that $L_i$ satisfies $L_i.C_i=0$ and $L_i.G>0$ for any pseudoeffective curve $[G] \in \overline{NE}(S)$ such that $[G] \not\in\mathbb{R}_{\geq0}[C_i]$. We identify $E$ with $S$. Let $\mathcal{L}_i$ be a pullback of $L_i$ by $\pi$ and $\mathcal{F}_i:= \phi^* (\mathcal{O}_{X_0}(1))\otimes \mathcal{L}_i$. We show that $\mathbb{R}_{\geq0}[C_i] \subseteq \overline{NE}(X)$ is an extremal ray. Since $(K_X+E)|_{E}\sim K_E$, it holds that $(K_X + E).C_i=-1$. By the cone theorem for dlt pairs, there exist finitely many $(K_X+E)$-negative extremal rays $R_k$ such that $[C_i]- [D] \in \sum R_k$ for some $[D] \in \overline{NE}(X)_{K_X+E=0}$. It holds that $\mathcal{F}_i .D =\mathcal{F}_i .R_k =0$ for all $k$ since $\mathcal{F}_i. C_i=0$ and $\mathcal{F}_i$ is a nef line bundle.
We see that, if an effective $1$-cycle $C$ on $X$ satisfies $\mathcal{F}_i .C=0$, then $C=\alpha C_i$ for some $\alpha \geq 0$ by the construction of $\mathcal{F}_i$. Thus, any generator of $R_k$ is equal to $\alpha_k C_i$ for some $\alpha_k \geq 0$. Hence $\mathbb{R}_{\geq0}[C_i] \subseteq \overline{NE}(X)$ is an extremal ray. It is clear to see that $\mathbb{R}_{\geq0}[C_i]\not=\mathbb{R}_{\geq0}[C_j]$. Thus the claim holds.
\end{proof}
\end{eg}

\end{document}